\documentclass[oneside,12pt]{amsart}
\usepackage{amssymb, amscd}
\usepackage[all]{xy}

\newtheorem{thm}{Theorem}[section]
\newtheorem{cor}[thm]{Corollary}
\newtheorem{lem}[thm]{Lemma}
\newtheorem{defn}[thm]{Definition}
\newtheorem{prop}[thm]{Proposition}
\newtheorem{examp}[thm]{Example}

\newcommand{\del}[2]{{}}

\title{Toroidal Orbifolds, Gerbes and Group Cohomology}

\author{Alejandro Adem }
\address{Department of Mathematics, University of Wisconsin, Madison WI 53706, USA}%
\email{adem@math.wisc.edu}%
\author{Jianzhong Pan }%
\address{Institute of Mathematics, Academia Sinica, Beijing 100080, China}%
\email{pjz@math03.math.ac.cn}%

\thanks{The first author was partially supported by
the NSF and the second author was partially supported by NSFC project
19701032 }%
\subjclass{}%
\keywords{orbifolds, gerbes, group cohomology}%

\date{\today}

\begin{document}
\begin{abstract}
In this paper we compute the integral cohomology of certain
semi--direct products of the form $\mathbb Z^n\rtimes G$, arising
from a linear $G$ action on the $n$--torus, where $G$ is a 
finite group. The main application
is the complete calculation of torsion gerbes for six dimensional
examples arising in string theory.
\end{abstract}

\maketitle
\tableofcontents
\section{Introduction}

During the last several years, there has been a surge of interest in
the mathematics involved in studying
various twisting processes in string theory.
One of the earliest ones was the notion of a $B$-field on a smooth
manifold. A $B$-field can be interpreted 
topologically as a class in
$H^3(X, \mathbb Z)$. In particular  
the torsion classes, or the
so--called flat B-fields, 
appear to
be geometrically significant.
Another source of motivation comes from 
\emph{stringy orbifold
theory}. The twisting process on orbifolds is much more subtle
due to the presence of a local finite group action. In the case
of a global quotient orbifold $X/G$, the first notion of
twisting--discrete torsion--was proposed by Vafa \cite{v}: a
discrete torsion is a cohomology class $\alpha \in H^2(G,
U(1))$, where $G$ is a finite group. 
It has also been
suggested that discrete torsion relates to desingularization. 
Recall that there are two
methods to remove singularities, resolution and deformation. Both
play important roles in the theory of Calabi-Yau 3-folds. A smooth
manifold $Y$ obtained from an orbifold $X$ via a sequence of
resolutions and deformations is called a \emph{desingularization}
of $X$. In string theory, we also require the resolutions to be
crepant resolutions. It is known that such a desingularization may
not exist in dimensions higher than three. In this case, we allow
our desingularization to be an orbifold. Vafa and Witten 
\cite{vw} proposed 
discrete torsion as a parameter for deformation. However, this
proposal immediately runs into trouble because the number of
desingularizations is much larger than the number of discrete
torsion elements. For example, D. Joyce \cite{joyce} constructed five
different desingularizations of an orbifold of the form
$\mathbb T^6/\mathbb{Z}/4$ (see Example \ref{ex1} below), while
$H^2(\mathbb{Z}/4, U(1))=0$. To count these ``missing''
desingularizations, Ruan \cite{inner} proposed the notion of an inner
local system. 
In the framework of gerbes, an inner local system can be
interpreted as
a holonomy line bundle (see \cite{lu}).

The computation of gerbes can be a difficult problem.
For the examples arising from orbifolds, 
almost all the known computations are from
discrete torsion. A priori we do not even know if gerbes capture
Joyce's example, which is what we would like to hope for.
It seems then that to understand gerbes better it will be 
necessary to
compute them for
physically
significant
examples. This is the
main motivation for our paper. 

Here we establish a 
cohomology computation
for groups arising as orbifold
fundamental groups associated to examples from geometry and string theory.
In particular our computations apply to six--dimensional examples arising
from finite group actions on tori, 
and we obtain some rather interesting 
examples of \textsl{torsion gerbes}, which should have
physical significance.
From the point of view of group cohomology, we are extending 
previous work (see \cite{adem}), where a number of calculations
for equivariant cohomology were explicitly obtained. 
Our mathematical interest was rekindled by more recent
contributions related to
physics, such as the work appearing in
\cite{bd}, where a comprehensive description of computations for
four--dimensional examples was provided. 

We now outline the basic mathematical framework for our paper. Let
$G\subset GL_n(\mathbb Z)$ denote a finite subgroup; then it will
act linearly on $\mathbb R^n$ preserving a lattice $\mathbb Z^n\subset
\mathbb R^n$. Hence there is an induced \textsl{linear} action on the
torus $\mathbb T^n$. Our basic goal will be to compute the equivariant
cohomology $H^*(EG\times_G \mathbb T^n,\mathbb Z)$. 
Now the Borel construction $EG\times_G \mathbb T^n$ will be aspherical, whence
it is an Eilenberg-MacLane space of type
$K(\pi_1(EG\times_G \mathbb T^n),1)$. If we let
$\Gamma = \pi_1(EG\times_G \mathbb T^n)$, then the problem we are raising is that
of computing the low dimensional cohomology of the discrete group
$\Gamma$. Note that $\Gamma$ can also be identified with the
\textsl{orbifold fundamental group of $\mathbb T^n\to \mathbb T^n/G$}.
As a consequence of our calculations, we will be able to compute
the group of isomorphism classes of (flat) gerbes for these orbifolds,
as these invariants can be identified with the (second)
third equivariant cohomology group of the manifold with ($U(1)$)
integral coefficients. In other words the invariants we will be
interested in computing are $FGb(\mathbb T^n/G)=
H^2(EG\times_G \mathbb T^n, U(1))$ and
$Gb(\mathbb T^n/G)=H^3(EG\times_G \mathbb T^n, \mathbb Z)$.

For all of these actions, $[0]\in\mathbb R^n/\mathbb Z^n$ is a fixed point.
This implies that there is a geometric section for $EG\times_G\mathbb T^n$,
hence also an algebraic 
section for $\Gamma\to G$. Therefore $\Gamma \cong \pi_1 (\mathbb T^n)\rtimes G
= \mathbb Z^n\rtimes G$, a semi--direct product. 
An important ingredient here will be the
detailed structure of the $\mathbb ZG$--lattice $\mathbb
Z^n$.\footnote{If $G$ is a finite group, a $\mathbb ZG$--lattice
is a $\mathbb ZG$--module which happens to be a free abelian group
of finite rank.} Indeed, the semi--direct product is completely
determined by this $\mathbb Z [G]$--module structure.

The following theorem is the basic computational device which
we will use:

\begin{thm}
If $M$ is a $\mathbb ZG$--lattice which restricts on each
$p$--Sylow subgroup to a direct sum of
finitely many one and two dimensional integral
representations, 
then for all integers $k\ge 0$, we have 
$$H^k(M\rtimes G,\mathbb Z)\cong
\bigoplus_{i+j=k}~H^i(G,H^j(M,\mathbb Z)).$$
\end{thm}

Having established this, we turn our attention to examples of 
interest in physics.
Computations for
four dimensional models were explicitly
described in \cite{bd}. 
However, for all of those examples, the odd dimensional
cohomology vanishes, hence no interesting
torsion gerbes can appear.
In contrast, the  six--dimensional examples
are far more interesting and much harder
to compute. In particular we consider the following important
examples studied in \cite{joyce} and \cite{vw}.

\begin{examp}\label{ex1}
Let $Y_1=T^6/\mathbb{Z}/4$ where 
the action is induced from that on $\mathbb{C}^3$,
given by
$$\tau(z_1,z_2,z_3)=(-z_1,iz_2,iz_3).$$
\end{examp}
\begin{examp}
Let $Y_2=T^6/{\mathbb{Z}/2}\times\mathbb Z/2$ where
the action
is induced from that on $\mathbb C^3$, given by
$$\tau_1(z_1,z_2,z_3)=(-z_1,-z_2, z_3), \text{  } 
\tau_2(z_1,z_2,z_3)=(-z_1,z_2,-z_3).$$
\end{examp}

Our techniques yield the following complete calculations for
the corresponding groups of gerbes:

$$FGb(Y_1)=U(1)^5 \oplus (\mathbb{Z}/2)^4 ,
\text{  }Gb(Y_1)=\mathbb{Z}^4 \oplus (\mathbb{Z}/2)^4$$

$$FGb(Y_2)=U(1)^3\oplus (\mathbb{Z}/2)^{18} , \text{
}Gb(Y_2)=\mathbb{Z}^8 \oplus (\mathbb{Z}/2)^{18}$$
Note the rather
mysterious torsion subgroup $(\mathbb{Z}/2)^4$ for
$Y_1$
and $(\mathbb{Z}/2)^{18}$ for
$Y_2$. We
expect that it will have physical significance, and it would
be interesting to
build a geometric model for these torsion
gerbes. Their role in twisted orbifold invariants 
(see \cite{inner}) also deserves
further study.
However we should point out that a
consequence of our computations is the conclusion
that it is too naive to expect that gerbes can capture the different 
desingularizations
of an orbifold. Instead we see that cohomology calculations for
the classifying space of an orbifold supply intricate invariants
which will measure the complicated nature of the possible resolutions, but
are far from determining them.

The methods used in this paper are from the cohomology of
groups. Applying a construction from \cite{br} we will manufacture
a very simple projective resolution for certain
semidirect products of the
form $\mathbb Z^n\rtimes G$, where $G$ is a finite group.
The key point will be that the resolution will allow us to
prove--at the cochain level--the
collapse (at the $E_2$ stage) of the Lyndon--Hochschild--Serre
spectral sequence for the group extension. The methods are 
actually very similar
to those used to compute the cohomology of wreath products (Nakaoka's
Theorem, see \cite{evens}, page 50). 

We should mention that computing the cohomology of groups such as these
(including crystallographic groups), is an important problem in the cohomology
of groups. 
More specifically, the exact
determination of the
Schur multiplier $H^2(\mathbb Z^n\rtimes G, U(1))$ seems like an
important basic
question.

The authors would like to thank D. Benson,
A. Klemm and Y.Ruan for their comments and
advice.

\section{Cohomology of Semidirect Products}

General references for group cohomology are the books by Adem--Milgram
\cite{am}, Brown \cite{brown}
and Evens \cite{evens}.
We will assume a basic familiarity with
the material on group extensions and their cohomology. 

The key input which we will use is a result due to 
Brady \cite{br}, which will allow us to construct an extremely
practical projective resolution and from there compute the cohomology
of the groups we will be interested in. The basic approach will be
the following: given certain conditions on
the $\mathbb ZG$--lattice $M$, the cohomology of the semidirect
product $M\rtimes G$ can be completely computed.\footnote{Recall
that given a split extension with a prescribed action on an
abelian kernel, the resulting group is uniquely determined up to
isomorphism.}

Let $\Gamma=M\rtimes G$ where $G$ acts on the group $M$ by a
homomorphism $\phi:G \to Aut(M)$. We write $m^g$ for $\phi(g)m$ where 
$m\in M , g\in G$.
\begin{defn}\label{T:comp}
Given a free resolution $\epsilon:F \to \mathbb{Z}$ of
$\mathbb{Z}$ over $\mathbb Z[M]$, 
we say that it admits an action of $G$
compatible with $\phi$ if for all $g\in G$ there is an
augmentation-preserving chain map $\tau(g):F \to F$ such that
\begin{enumerate}
 \item{$\tau(g)[m\cdot f]=m^g\cdot [\tau(g)f]$
for all $m\in M$  and $f\in F$}, 
\item{$\tau(g)\tau(g')=\tau(gg')$ for all $g,g'\in G$},
\item{$\tau (1) = 1_F$.}
\end{enumerate}
\end{defn}

If such an action exists we can give $F$ a $\Gamma$-module structure as
follows. If $\gamma\in \Gamma$ then $\gamma$ can be expressed uniquely as 
$\gamma=mg$,
with $m\in M$ and $g\in G$. We set $\gamma\cdot f=(mg)\cdot f=m\cdot
\tau(g)f$. Note that given any $G$--module $L$, this inflates to 
a $\Gamma$--action on $L$ via the projection $\Gamma\to G$.

\begin{prop}\label{T:exis}
If $\epsilon:F \to \mathbb{Z}$ is a free resolution of
$\mathbb{Z}$ over $\mathbb Z [M]$, which admits an action of $G$ compatible
with $\phi$, and $\epsilon':P \to \mathbb{Z}$ is any free
resolution of $\mathbb{Z}$ over $G$, then $\epsilon' \otimes
\epsilon : P \otimes F \to \mathbb{Z}$ is a free resolution of
$\mathbb{Z}$ over  $\mathbb Z\Gamma$.
\end{prop}
\begin{proof}
It suffices to show that each chain module is $\mathbb Z\Gamma$--free.
Indeed if $F_i$ restricts on $M$ to a free module of rank
$n_i$, then we have 
$$\mathbb Z G\otimes F_i\cong \mathbb Z [\Gamma /M]\otimes
F_i\cong \mathbb Z\Gamma\otimes_M (\mathbb Z [M])^{n_i}
\cong (\mathbb Z \Gamma)^{n_i}.$$
\end{proof}

In the case $\Gamma=\mathbb{Z}^n \rtimes G$, a particular
free resolution for
$\mathbb{Z}^n $ can be chosen as follows. Let
$\partial: \mathbb Z [\mathbb Z]\to \mathbb Z[\mathbb Z]$
be given by $\partial (u) = (1-x)u$, where $x$ is a generator for
the infinite cyclic group. Then one can verify that $\partial$ is
injective, with cokernel the trivial module $\mathbb Z$. Denote
this free resolution by $D$; take $n$ copies of it, say denoted by
$D_1, D_2, \dots , D_n$. Now let $F=D_1\otimes D_2\otimes\dots
\otimes D_n$; this will be a free resolution for the trivial
module $\mathbb Z$ over the group ring $\mathbb Z [\mathbb Z^n]$.
We observe the following: $Hom_{\mathbb Z [\mathbb Z]}(D,\mathbb
Z)$ has trivial differentials. Using the Kunneth formula, we infer
that the analogous fact
will hold for the resolution $F$. Hence we have that
the cochain complex for computing the cohomology $H^*(\mathbb
Z^n,\mathbb Z)$ is in fact additively isomorphic to the
cohomology. We shall call this a \textbf{special resolution}.

Now suppose we are given $\Gamma=\mathbb Z^n\rtimes G$ and the
existence of a compatible action of $G$ on this special resolution
$F$ for $\mathbb Z^n$. Then we can apply the 
ideas\footnote{Our proof is
very similar to Evens' proof of Nakaoka's Theorem on the cohomology
of wreath products (see \cite{evens}, pages 19 and 50). We are grateful
to Dave Benson for this observation.} 
above as
follows.

\begin{thm}\label{T:2.3}
Let $\epsilon:F \to \mathbb{Z}$ be a special free resolution of
$\mathbb{Z}$ over $\mathbb{Z}^n $ and suppose that there is a
compatible action of $G$ on $F$. Then for all integers $k\ge 0$, 
we have
$$H^k(\mathbb{Z}^n \rtimes
G,\mathbb Z)=\bigoplus_{i+j=k}H^i(G,H^j(\mathbb{Z}^n,\mathbb
Z)).$$
\end{thm}
\begin{proof}
Let $\epsilon' \otimes \epsilon : P \otimes F \to \mathbb{Z}$ be
the free resolution of $\mathbb{Z}$ over $\Gamma$ which is given
by Proposition \ref{T:exis}. Then the cohomology of $\Gamma$ can
be computed as the cohomology of $Hom_{\Gamma}(P\otimes F, \mathbb
Z)$. Note however that the normal subgroup $\Gamma'=\mathbb Z^n$
acts trivially on $P$, whence we obtain the isomorphism
$Hom_{\Gamma}(P\otimes F, \mathbb Z)\cong
Hom_{\Gamma/\Gamma'}(P, Hom_{\Gamma'}(F,\mathbb Z))$. 
Hence our
computation has reduced to calculating the $G$
equivariant cohomology of a $\mathbb ZG$--cochain complex with trivial
differentials, from which the result immediately follows.
\end{proof}

To conclude the section, we will state three simple but useful
lemmas which facilitate the construction of compatible actions;
their proofs are straightforward.

\begin{lem}\label{T:split}
 If $\epsilon_i:F_i \to \mathbb{Z}$ is a projective
 $M_i$-resolution of
 $\mathbb{Z}$ for $i=1,2$, then $\epsilon_1
 \otimes \epsilon_2:F_1 \otimes F_2 \to
 \mathbb{Z}$ is a projective $M$-resolution of $\mathbb{Z}$,
 where $M=M_1 \times
 M_2$. If, furthermore, a group $G$ acts compatibly on $F_i$ by $\tau_i$ for
 $i=1,2$, then a compatible action of $G$ on $\epsilon_1 \otimes \epsilon_2:F_1 \otimes F_2 \to
 \mathbb{Z}$ is given by $\tau(g)(f_1 \otimes f_2)=
 \tau_1(g)(f_1) \otimes \tau_2(g)(f_2)$
\end{lem}

\begin{lem}\label{T:1.5}
 If $G=G_1 \times G_2$, $M$ is a $G_i$-module for $i=1,2$, and $G_i$ acts
 compatibly by $\tau_i$
 on a $M$-resolution of
 $\mathbb{Z}$,
 $\epsilon:F \to \mathbb{Z}$ such that $\tau_1(g_1)\tau_2(g_2)=\tau_2(g_2)\tau_1(g_1)$, then
 a compatible action of $G$ on $\epsilon:F \to
 \mathbb{Z}$ is given by $\tau(g)(f)=\tau_1(g)  \tau_2(g)f$
\end{lem}

\begin{lem}\label{T:1.6}
 If $M$ is a $G_1$-module, $\pi:G_2 \to G_1$ a group homomorphism, and
 $\epsilon:F \to \mathbb{Z}$ is a $M$-resolution of
 $\mathbb{Z}$ such that $G_1$ acts compatibly on it by $\tau'$, then
 $G_2$ also acts compatibly on it via
 $\tau(g)f=\tau'(\pi(g))f$ for any $g\in G_2$
\end{lem}

\section{Constructing Compatible Actions}
Consider the special $\mathbb Z [\mathbb{Z}]$-resolution of
$\mathbb{Z}$ :

\[0 \to \mathbb{Z}[\mathbb{Z}]a_1 \to \mathbb{Z}[\mathbb{Z}]a_0 \to
\mathbb{Z}\]
with generators $a_i$ in degree $i$ for $i=0,1$.
If we tensor it with another copy of the same resolution, but now with
generators $b_0, b_1$, we obtain the explicit special resolution for
the group $\mathbb Z^2= \mathbb Z ~x_1\oplus \mathbb Z ~x_2$:

$$
\mathbb Z[\mathbb Z^2]a_1\otimes b_1\to
\mathbb Z[\mathbb Z^2]a_0\otimes b_1\oplus
\mathbb Z[\mathbb Z^2]a_1\otimes b_0\to
\mathbb Z[\mathbb Z^2]a_0\otimes b_0
$$
where the boundary maps are given by
$$\partial (a_1\otimes b_1) = (1-x_1)a_0\otimes b_1 - (1-x_2)a_1\otimes b_0$$
$$\partial (a_0\otimes b_1) = (1-x_2)a_0\otimes b_0,
~~~~\partial (a_1\otimes b_0) = (1-x_1)a_0\otimes b_0.$$

In what follows we will be constructing compatible finite group
actions on this explicit complex. The action on $\langle x_1,
x_2\rangle$ is given as follows. If $A\in GL_2(\mathbb Z)$, then
we identify $x_1$ with $(1,0)$ and $x_2$ with $(0,1)$
respectively. Then the action is of the form
$$(1,0)\mapsto
(1,0)A=(a_{11}, a_{12}), \;\;\;\;\;(0,1)\mapsto (0,1)A =(a_{21},
a_{22})$$
which multiplicatively becomes
$$x_1\mapsto
x_1^{a_{11}}x_2^{a_{12}},\;\;\;\;x_2\mapsto
x_1^{a_{21}}x_2^{a_{22}}.$$

Given a generator $t\in G$, we can identify it with a $2\times 2$ matrix
$A\in GL_2(\mathbb Z)$ as above. We now describe an explicit procedure
for constructing $\tau (t):F\to F$. First we observe that as we want it to
be augmentation--preserving, we assume that the generator $a_0\otimes b_0$ 
is fixed.
Next we set up the algebraic equations which must be solved to define
$\tau (t)$ as a chain map.

For brevity we write $\tau$ instead of $\tau (t)$. We have

$$\tau (a_0\otimes b_1) = q_{01} a_0\otimes b_1 + q_{10} a_1\otimes b_0,
\;\;\;
\tau (a_1\otimes b_0) = r_{01} a_0\otimes b_1 + r_{10} a_1\otimes b_0$$
for elements $r_{ij}, q_{ij}\in \mathbb Z [\mathbb Z^2]$. Now the equation
$\partial \tau = \tau\partial$ is equivalent to the following algebraic conditions 
on these coefficients:

$$q_{01}(1-x_2) + q_{10}(1-x_1) = 1-x_1^{a_{21}}x_2^{a_{22}},\;\;\;
r_{01}(1-x_2) + r_{10}(1-x_1) = 1 - x_1^{a_{11}}x_2^{a_{12}}.$$
Similarly let's assume that $\tau (a_1\otimes b_1) = q_{11}a_1\otimes b_1$.
Then the chain complex condition implies the following equations:

$$q_{11}(1-x_1) = (1-x_1^{a_{11}}x_2^{a_{12}})q_{01}- 
		  (1-x_1^{a_{21}}x_2^{a_{22}})r_{01}$$
$$-q_{11} (1-x_2) = (1-x_1^{a_{11}}x_2^{a_{12}})q_{10} -
		    (1-x_1^{a_{21}}x_2^{a_{22}})r_{10}.$$

\bigskip

\begin{thm}\label{T:3.1}
Let $M$ be a $\mathbb ZG$--lattice which is a direct sum of
indecomposable $G$-submodules of $\mathbb{Z}$-rank at most two.
Then there is a compatible action of $G$ on the special
$M$-resolution of $\mathbb{Z}$ if $G$ is as follows:
\begin{enumerate}
\item $G=\mathbb{Z}/n$ for $n=2,3,4$,

\item $G=\mathbb{Z}/2 \oplus \mathbb{Z}/2$,

\item $G=D_8$, the dihedral group of order $8$.

\end{enumerate}
\end{thm}
\begin{proof}
We will be considering isomorphism classes of $\mathbb
ZG$--lattices. Note that by Lemma \ref{T:split} it suffices to
prove the theorem when the module is indecomposable of rank one or
two. On the other hand, for rank one modules, it suffices to prove
the theorem for $G=\mathbb{Z}/2$ by Lemma \ref{T:1.6}. Now we
proceed to prove the theorem case by case.

\noindent (1)~$G=\mathbb{Z}/2$ with generator $t$;
the only nontrivial indecomposable $G$-modules of $\mathbb{Z}$-rank
at most $2$ are given by the matrix representing $t$: $(-1)$ and
\[
A= \left( \begin{array}{ll} 0 & 1 \\ 1 & 0 \end{array} \right).
\]
For the first case, $\tau(t)(a_0)=a_0$ and
$\tau(t)(a_1)=-x^{-1}a_1$. It is easy to verify that this gives a
compatible action of $\mathbb{Z}/2$ on the above resolution.

For the second case, we use the special $\mathbb Z[\mathbb{Z}^2]$-resolution
of $\mathbb{Z}$ described above, with free generators denoted by
$a_i\otimes b_j$ where $i,j=0,1$
For this example,
the equations above admit the following simple solutions:
$q_{01}=0$, $q_{10}=1$, $r_{01}=1$, $r_{10}=0$, $q_{11}=-1$.
Hence the following definition defines a chain map satisfying
the first compatibility condition:
$$\tau(t)(a_0 \otimes b_0)=a_0 \otimes b_0, \;\;\;\; 
\tau(t)(a_1 \otimes b_0)=a_0 \otimes b_1 $$
$$\tau(t)(a_0 \otimes b_1)=a_1 \otimes b_0, \;\;\;\;
\tau(t)(a_1 \otimes b_1)=-a_1 \otimes b_1$$
Moreover in this case one can easily verify that this
gives a compatible action of $\mathbb{Z}/2$ on the above
resolution, i.e., $\tau(t)^2=Id:F \to F$.

\noindent (2)~$G=\mathbb{Z}/3$ with
generator $t$;
the only nontrivial indecomposable
$\mathbb ZG$-module of $\mathbb{Z}$-rank two can be
represented by the
matrix 
$$A= \left( \begin{array}{ll} 0 & -1 \\ 1 & -1 \end{array} \right).$$
In this case the solutions to the equations are incorporated
into the following formula for $\tau$:
$$
\begin{array}{ccc}
\tau(t)(a_0 \otimes b_0)&=&a_0 \otimes b_0 \\
\tau(t)(a_1 \otimes b_0)&=&-x_2^{-1}a_0 \otimes b_1\\
\tau(t)(a_0 \otimes b_1)&=&a_1 \otimes b_0-x_1x_2^{-1}a_0 \otimes b_1\\
\tau(t)(a_1 \otimes b_1)&=&x_2^{-1}a_1 \otimes b_1.
\end{array}
$$
We will give a complete proof that this indeed gives a compatible
action of $\mathbb{Z}/3$ on the special resolution.
We need to check condition 2) in definition \ref{T:comp} which
in this case is $\tau(t)^3=Id$.

$$\begin{array}{ccc}
\tau(t)^3(a_1 \otimes b_0)&=&\tau(t)^2(-x_2^{-1}a_0 \otimes
b_1)\\
&=&-\tau(t)(x_2x_1^{-1}(a_1 \otimes b_0-x_1x_2^{-1}a_0 \otimes
b_1))\\
&=&\tau(t)(a_0 \otimes b_1-x_2x_1^{-1}a_1 \otimes b_0)\\
&=&a_1 \otimes b_0-x_1x_2^{-1}a_0 \otimes b_1+x_1x_2^{-1}a_0
\otimes b_1\\
&=&a_1 \otimes b_0\\
&&\\
\tau(t)^3(a_0 \otimes b_1)&=&\tau(t)^2(a_1 \otimes
b_0-x_1x_2^{-1}a_0 \otimes b_1)\\
&=&\tau(t)(-x_2^{-1}a_0 \otimes b_1-x_1^{-1}(a_1 \otimes
b_0-x_1x_2^{-1}a_0 \otimes b_1))\\
&=&\tau(t)(-x_1^{-1}a_1 \otimes
b_0)\\
&=&-x_2(-x_2^{-1}a_0 \otimes b_1)=a_0 \otimes b_1\\
&&\\
\tau(t)^3(a_1 \otimes b_1)&=&\tau^2(x_2^{-1}a_1 \otimes
b_1)\\
&=&\tau(t)(x_2x_1^{-1}x_2^{-1}a_1 \otimes b_1)\\
&=&\tau(t)(x_1^{-1}a_1 \otimes b_1)\\
&=&x_2x_2^{-1}a_1 \otimes b_1)\\
&=&a_1 \otimes b_1
\end{array}$$

\noindent (3)~$G=\mathbb{Z}/4$ with generator $t$;
the only faithful indecomposable $\mathbb{Z}/4$--lattice of
$\mathbb{Z}$-rank $2$ is given by the matrix representing $t$:

\[
A= \left( \begin{array}{ll} 0 & 1 \\ -1 & 0 \end{array} \right)
\]

Again in this case we can solve the equations and
so we have:
$$
\begin{array}{ccc}\tau(t)(a_0 \otimes b_0)&=&a_0 \otimes b_0\\
\tau(t)(a_1 \otimes b_0)&=&a_0 \otimes b_1 \\
\tau(t)(a_0 \otimes b_1)&=&-x_1^{-1}a_1 \otimes b_0\\
\tau(t)(a_1 \otimes b_1)&=&x_1^{-1}a_1 \otimes b_1
\end{array}
$$
We leave the verification that $\tau (t)^4=1$ to the reader.

\noindent (4)~$G=\mathbb{Z}/2 \oplus \mathbb{Z}/2$ with generators $t_1,t_2$:
the only faithful indecomposable $\mathbb ZG$--lattice of
$\mathbb{Z}$-rank $2$ is given by the matrices representing
$t_1,t_2$:

$$
 A_1= \left( \begin{array}{ll} 0 & 1 \\ 1 & 0 \end{array}
\right), \;\;A_2= \left( \begin{array}{ll} -1 & 0\\ 0 & -1
\end{array} \right). $$
In this case, by part (1) above, and Lemma \ref{T:split},
$\tau(t)$ can be defined as follows:

$$
\begin{array}{ccc}
\tau(t_2)(a_1 \otimes b_0)&=&-x_1^{-1}a_1 \otimes b_0\\

\tau(t_2)(a_0 \otimes b_1)&=&-x_2^{-1}a_0 \otimes b_1\\

\tau(t_1)(a_1 \otimes b_0)&=&a_0 \otimes b_1\\

\tau(t_1)(a_0 \otimes b_1)&=&a_1 \otimes b_0\\

\tau(t_2)(a_1 \otimes b_1)&=&x_1^{-1}x_2^{-1}a_1 \otimes b_1\\

\tau(t_1)(a_1 \otimes b_1)&=&-a_1 \otimes b_1.
\end{array}$$

By Lemma \ref{T:1.5}, it suffices to prove that
$\tau(t_1)\tau(t_2)=\tau(t_2)\tau(t_1)$

$$\tau(t_1)\tau(t_2)(a_1 \otimes b_0)=\tau(t_1)(-x_1^{-1}a_1
\otimes b_0) =-x_2^{-1}a_0 \otimes b_1$$

$$\tau(t_2)\tau(t_1)(a_1 \otimes b_0)=\tau(t_2)(a_0 \otimes
b_1)=-x_2^{-1}a_0 \otimes b_1$$

Similarly we see that $\tau(t_1)\tau(t_2)(a_0 \otimes
b_1)=\tau(t_2)\tau(t_1)(a_0 \otimes b_1)$. Next we check

$$\tau(t_1)\tau(t_2)(a_1 \otimes
b_1)=\tau(t_1)(x_1^{-1}x_2^{-1}a_1 \otimes b_1)
=-x_1^{-1}x_2^{-1}a_1 \otimes b_1$$

$$\tau(t_2)\tau(t_1)(a_1 \otimes b_1)=\tau(t_2)(-a_1 \otimes
b_1) =-x_1^{-1}x_2^{-1}a_1 \otimes b_1.$$

\noindent (5)~Dihedral group case 
$G=<t_1,t_2| t_1^4=t_2^2=e,t_2t_1t_2=t_1^{-1}>$;
the only faithful indecomposable $\mathbb ZG$--lattices of
$\mathbb{Z}$-rank $2$ are given by the matrices representing
$t_1,t_2$:

$$ A_1= \left( \begin{array}{ll} 0 & 1 \\ -1 & 0 \end{array}
\right) ,\;\;A_2= \left( \begin{array}{ll} 1 & 0\\ 0 & -1
\end{array} \right) $$

$$ B_1= \left( \begin{array}{ll} 0 & 1 \\ -1 & 0 \end{array}
\right),\;\;B_2= \left( \begin{array}{ll} 0 & 1\\ 1 & 0
\end{array} \right) $$
In the first case we solve the equations to obtain
$$\begin{array}{ccc}
\tau(t_1)(a_1 \otimes b_0)&=&a_0 \otimes b_1\\

\tau(t_1)(a_0 \otimes b_1)&=&-x_1^{-1}a_1 \otimes b_0\\

\tau(t_2)(a_1 \otimes b_0)&=&a_1 \otimes b_0\\

\tau(t_2)(a_0 \otimes b_1)&=&-x_2^{-1}a_0 \otimes b_1\\

\tau(t_1)(a_1 \otimes b_1)&=&x_1^{-1}a_1 \otimes b_1\\

\tau(t_2)(a_1 \otimes b_1)&=&-x_2^{-1}a_1 \otimes b_1
\end{array}$$

In the second case we have:

$$\begin{array}{ccc}
\tau(t_1)(a_1 \otimes b_0)&=&a_0 \otimes b_1\\

\tau(t_1)(a_0 \otimes b_1)&=&-x_1^{-1}a_1 \otimes b_0\\

\tau(t_2)(a_1 \otimes b_0)&=&a_0 \otimes b_1\\

\tau(t_2)(a_0 \otimes b_1)&=&a_1 \otimes b_0\\

\tau(t_1)(a_1 \otimes b_1)&=&x_1^{-1}a_1 \otimes b_1\\

\tau(t_2)(a_1 \otimes b_1)&=&-a_1 \otimes b_1
\end{array}$$

\end{proof}

We will use restriction to the $p$--Sylow subgroups
to obtain a cohomology calculation.

\begin{prop}
Let $\Gamma = M\rtimes G$ denote an extension of a finite
group by a finitely generated free abelian group. If the 
Lyndon--Hochschild--Serre spectral sequence restricted to every
$p$--Sylow subgroup $G_p\subset G$ collapses at $E_2$, without
extension problems, then the
same is true for the LHS spectral sequence of the original extension.
\end{prop}
\begin{proof}
The restriction induces an
embedding $H^i(G,H^j(M))_{(p)}\subset 
H^i(G_p, H^j(M))$ for all $i>0$, whereas $H^j(M)^G\subset 
H^j(M)^{G_p}$. Hence we deduce that the differentials are
trivial at $p$, for all primes $p$ dividing $|G|$. As the images always
lie in torsion modules annihilated by $|G|$, we deduce that 
there can be no differentials. Similarly if there are no extension
problems p--locally, they cannot occur integrally.
\end{proof}

\begin{cor}\label{LHS}
If $M\rtimes G$ admits a compatible action for all $p$--Sylow
subgroups $G_p\subset G$, then the LHS spectral sequence for the
extension collapses at $E_2$ without extension problems.
\end{cor}

We can now prove our main result.
\begin{thm}\label{mainthm}
If $M$ is a $\mathbb ZG$--lattice which 
restricts on each $p$--Sylow subgroup of $G$ to
a direct sum 
of finitely many one and two dimensional integral
representations, 
then for all integers $k\ge 0$ we have
$$H^k(M\rtimes G,\mathbb Z)\cong
\bigoplus_{i+j=k}~H^i(G,H^j(M,\mathbb Z)).$$
\end{thm}
\begin{proof}
We need only consider $p$--groups and
integral representations of rank one
or two which are faithful.
It is clear that $GL_1(\mathbb Z)\cong\mathbb Z/2$, whereas the complete
list of finite $p$--subgroups of $GL_2(\mathbb Z)$ is given by

\begin{itemize}

\item The cyclic subgroups $\mathbb Z/2$, $\mathbb Z/3$,
$\mathbb Z/4$,

\item $\mathbb Z/2\times\mathbb Z/2$,

\item The dihedral group $D_8$ of order eight.

\end{itemize}
We have constructed compatible actions for all of the
faithful integral representations of these groups having rank
one or two.
Hence by Corollary \ref{LHS}, we obtain the desired
decomposition theorem.
\end{proof}

\section{Toroidal Orbifolds}

In this section we will outline how integral representations 
of finite groups naturally give rise to group actions on tori.
The resulting orbifolds appear as examples of interest in
physics and we will provide explicit computations in the next
section.

Assume that $G$ is a finite group and that we are given a faithful
integral representation $G\subset GL_n(\mathbb Z)$. This will induce
a smooth action of $G$ on $\mathbb R^n$ via matrices. Given that
these matrices have integer entries, this will induce an action on
the quotient $X=\mathbb R^n/\mathbb Z^n$, where $\mathbb Z^n\subset
\mathbb R^n$ is the canonical lattice of rank $n$. Note that
the action will
always have a fixed point (indeed the zero vector is invariant).
The space $Y=X/G$ is what is often denoted a \textsl{global quotient},
and $X\to X/G$ carries a natural orbifold structure.

Let $EG$ denote a contractible free $G$--space, and consider the
Borel construction $EG\times_G X$. Projection on the first
coordinate yields a fibration
$X\to EG\times_G X\to BG$, where $BG=EG/G$ is the classifying
space of $G$. Note that both $X$ and $BG$ have contractible universal
covers, hence $EG\times_GX$ is aspherical. This space can be
regarded as the \textsl{classifying space of the orbifold}
and $\Gamma=
\pi_1 (EG\times_GX)$ is known as the \textsl{orbifold fundamental
group}. Note that the cohomology of this group is by definition the
cohomology of its classifying space, which in this case is precisely
$EG\times_GX$. 

The group $\Gamma$  
fits into a group extension of the form
$1\to M\to\Gamma\to G\to 1$
where $M\cong \mathbb Z^n$ has the $G$ action induced from the
representation $G\subset GL_n(\mathbb Z)$ i.e. $M$ is a $\mathbb Z[G]$
lattice. This extension arises from applying $\pi_1$ to the fibration
above. Note that the projection $EG\times_GX\to BG$ has a section, 
due to the fact that the $G$--action has a fixed point. Hence the extension
splits, and so $\Gamma\cong M\rtimes G$, \textsl{the semidirect
product} determined by the action of $G$ on $M$ (up to isomorphism).
We will be interested in knowing how many different splittings
this extension has. This is measured by the cohomology group
$H^1(G, M)$. It can also be identified with $M$--conjugacy
classes of subgroups $\tilde{G}\subset \Gamma$
mapping isomorphically onto
$G$. 

We now express this information geometrically. The $G$--action
on $X$ lifts to a $\Gamma$ action on $Z=\mathbb R^n$. As pointed out in
\cite{brown}, page 267 we have, for any $Q\subset G$:
$$X^Q = \bigsqcup_{H^1(Q, M)} 
Z^{\tilde{Q}}/M\cap N(\tilde{Q})$$
where $N(\tilde{Q})$ is the normalizer of $\tilde{Q}$.
We can also consider the coefficient sequence
$$0\to M\to M\otimes\mathbb R\to X\to 0$$
which gives rise to the exact sequence
$$0\to M^Q\to (M\otimes\mathbb R)^Q\to X^Q\to H^1(Q,M)\to 0.$$
From this we see that for every $Q\subset G$,
the $Q$--fixed point set in $X$ will always decompose
into a disjoint union of tori. 
Comparing the two expressions 
we see that $Z^{\tilde{Q}}$ is isomorphic to $\mathbb R^{n_Q}$, where
$n_Q$ is the rank of the free abelian group $M^Q$, the subgroup of $Q$--invariants,
and that $M\cap N(\tilde{Q})\cong \mathbb Z^{n_Q}$. Indeed we 
obtain the
formula
$$X^Q = \bigsqcup_{H^1(Q,M)}~(\mathbb S^1)^{rk_{\mathbb Z}[M^Q]}.$$
Note in particular that if $M^Q=\{0\}$, 
then $X^Q$ will be a disjoint collection of points, indexed
by $H^1(Q, M)$.

In differential geometry examples arise using complex coordinates.
More precisely if say $G\subset GL_m(\mathbb C)$ and there exists
a rank $2m$ lattice $\Lambda\subset \mathbb C^m$ such that $G\Lambda
\subset \Lambda$, then this will induce a $G$--action on 
$(\mathbb S^1)^{2m} = \mathbb C^m/\Lambda$. The action on the lattice
$\Lambda$ (assumed to be faithful) provides an embedding
$G\subset GL(\Lambda)$ from which we derive the same data as before
and can therefore construct a topologically equivalent model.

\section{Key Examples and Calculations}

In this section we will apply our techniques to certain very
specific examples arising from geometry and physics. 
Given an orbifold $X\to X/G$ our main object of interest
will be the group of gerbes associated to it.
We shall not dwell here on the definition of gerbes, but
use the following identification for the isomorphism
classes of gerbes and flat gerbes respectively:

$$Gb(X/G)=H^3(EG\times_G X, \mathbb Z),\;\;\;
FGb(X/G) = H^2(EG\times_G X,\mathbb U(1)).$$
For a finitely generated abelian group $A$, let
$A =T(A) \oplus F(A)$ denote the decomposition into torsion and torsion--free
components respectively. Then,
by the universal coefficient theorem, we have
$$Gb(X/G) = T(H_2(EG\times_G X,\mathbb Z))\oplus 
F(H_3(EG\times_G X,\mathbb Z))$$
$$FGb(X/G) = T(H_2(EG\times_G X,\mathbb Z))\oplus 
F(H_2(EG\times_G X,\mathbb Z))\otimes U(1)$$ 

For our examples, these calculations will amount to
computing the low dimensional cohomology of certain
discrete groups. In particular we will be computing
the \textsl{Schur multiplier} $H^2(M\rtimes G, \mathbb S^1)$.
This is a problem of independent interest in group cohomology,
and our methods indicate a practical strategy for computing
this invariant for many examples. 

First we introduce some general notation. If $M$ is a 
$\mathbb ZG$--lattice building the semidirect product
$M\rtimes G$, then the cohomology of the associated torus
can be described as an exterior algebra $\wedge^*(M^*)$
where $M^*$ is the dual of $M$. The action on these exterior
powers is determined by the original action on $M$. From the point
of view of group cohomology we should observe that for $G$ having
cyclic $p$--sylow subgroups, any $\mathbb ZG$--lattice will be
cohomologous to its dual. 

As would be expected, the most interesting examples, and the real motivation for 
this work, comes from calculations for six--dimensional orbifolds, where the
usual spectral sequence techniques become rather complicated. Here our
methods provide an important new ingredient that will allow us to compute
rigorously beyond the known range.
An important class of examples in physics arise 
from actions of a cyclic group
$\mathbb Z/N$ on $\mathbb T^6$. 
In our scheme, these come from
six-dimensional integral representations of $\mathbb Z/N$. However,
the constraints
from physics impose certain restrictions on them (see \cite{EK}).
If $\theta\in GL_6(\mathbb Z)$ is an element of order $N$, then
it can be diagonalized over the complex numbers.
The associated eigenvalues, denoted
$\alpha_1,\alpha_2,\alpha_3$, should satisfy
$\alpha_1\alpha_2\alpha_3 =1$, and in addition all of the
$\alpha_i\ne 1$. The first condition implies that
the orbifold $\mathbb T^6\to \mathbb T^6/G$ is a 
\textsl{Calabi--Yau orbifold}, and so admits 
a \textsl{crepant resolution}.
These more restricted representations
have been classified\footnote{In the language of
physics, they show that there exist 18 inequivalent $N=1$ supersymmetric
string theories on symmetric orbifolds of $(2,2)$--type without
discrete background.} in \cite{EK},
where it is shown that there are precisely 18 inequivalent lattices
of this type. Six of them
decompose as direct sums of sublattices of rank two or less, and so
we can apply Theorem \ref{mainthm} to obtain a computation of
the cohomology of the corresponding discrete groups. It would seem
likely that an extension of our methods should also apply to the other
examples on this list. We will provide a complete computation for
the following celebrated example, which appears in \cite{vw}.

\subsection{The case $Y_1=T^6/\mathbb{Z}/4$}
The action of $\mathbb{Z}/4$ on $T^6$ is induced from the action
of $\mathbb{Z}/4$ on $\mathbb{C}^3$, given by 
$$\kappa(z_1,z_2,z_3)=(-z_1,iz_2,iz_3).$$
The induced action of $\mathbb{Z}/4$  on $\mathbb{Z}^6$ is given
by the following decomposition of $\mathbb Z[\mathbb{Z}/4]$-lattices:
$$M=(M_1)^2 \oplus (M_2)^2 $$
where $M_1$ has rank one and $\kappa$ acts by  multiplication by
$-1$; and $M_2$ has rank two and $\kappa$ acts by the following rule
$\kappa(x)=-y ,\kappa(y)=x$, i.e. the standard embedding
$\mathbb Z/4\subset GL_2(\mathbb Z)$. We observe that this module
is a sum of lattices for which our methods apply, hence computing
the cohomology of the associated semi--direct product reduces to
determining the $E_2$ term of the LHS spectral sequence.

First note that $M\cong M^*$, so what we need to do is
determine the exterior powers of $M$. Note that
$\wedge^i(M)\cong \wedge^{6-i}(M)^*\cong \wedge^{6-i}(M)$.
Hence we only need to worry about $\wedge^2$ and $\wedge^3$.
We will use the distributive formula for exterior powers of
a direct sum. Now it is easy to see that $M_1\otimes M_1\cong
\mathbb Z$, the trivial module, whereas $M_1\otimes M_2\cong M_2$.
For the tensor product $M_2\otimes M_2$
we have a basis with action given by $x\otimes x\to y\otimes y$,
$x\otimes y\mapsto -y\otimes x$, $y\otimes x\mapsto -x\otimes y$ and
$y\otimes y\mapsto x\otimes x$. Hence $M_2\otimes M_2\cong
P\oplus P$, where $P\cong \mathbb Z[\mathbb Z/4]
\otimes_{\mathbb Z/2}\mathbb Z$
is the inflation of the regular representation
of $\mathbb Z/2$. 

Now we may compute: 
$$\wedge^2(M)\cong (\mathbb Z)^3\oplus (M_2)^4\oplus (P\oplus P),\;\;\;
\wedge^3(M)\cong (M_1)^4\oplus (M_2)^4\oplus (P)^4$$
Summarizing, we have:

\[H^i(\mathbb{Z}^6,\mathbb{Z})= \left\{ \begin{array}{ll}
\mathbb{Z} & \mbox{ if $i=0,6$ }\\ M_1^2 \oplus M_2^2 &
\mbox{ if $i=1,5$ } \\
\mathbb{Z}^3 \oplus (M_2)^4 \oplus (P)^2
& \mbox{ if $i=2,4$ }\\
(M_1)^4 \oplus (M_2)^4\oplus (P)^4
& \mbox{ if $i=3$ } \\ 0 &
\mbox{ if $i>6$ }
\end{array} \right.
 \]
Next we compute the terms $H^s(\mathbb Z/4, H^t (\mathbb Z^6, \mathbb Z))$.
These can be reduced 
to computing the cohomology of $\mathbb Z/4$ with coefficients
in the modules $\mathbb Z$, $M_1$, $M_2$ and $P$. By an elementary
calculation we have that

\[
H^i(\mathbb Z/4, \mathbb Z) = \left\{\begin{array}{ll}
\mathbb Z & \mbox{ if $i=0$ }\\0 & \mbox{ if $i$ is odd }\\
\mathbb Z/4 & \mbox{ if $i>0$ is even }\end{array}\right.
\]

\[
H^i(\mathbb Z/4, M_1) = \left\{\begin{array}{ll}
0 & \mbox{ if $i=0$ }\\\mathbb Z/2 & \mbox{ if $i$ is odd }\\
0 & \mbox{ if $i>0$ is even }\end{array}\right.
\]

\[
H^i(\mathbb Z/4, M_2) = \left\{\begin{array}{ll}
0 & \mbox{ if $i=0$ }\\\mathbb Z/2 & \mbox{ if $i$ is odd }\\
0 & \mbox{ if $i>0$ is even }\end{array}\right.
\]

\[
H^i(\mathbb Z/4, P) = \left\{\begin{array}{ll}
\mathbb Z & \mbox{ if $i=0$ }\\0 & \mbox{ if $i$ is odd }\\
\mathbb Z/2 & \mbox{ if $i>0$ is even }\end{array}\right.
\]

From this the cohomology can easily be written down:

\[H^i(\mathbb{Z}^6\rtimes \mathbb{Z}/4,\mathbb{Z})= \left\{ \begin{array}{ll}
\mathbb{Z} & \mbox{ if $i=0$ }\\ 0 &
\mbox{ if $i=1$ } \\
\mathbb{Z}^5 \oplus\mathbb{Z}/4\oplus (\mathbb Z/2)^4 & \mbox{ if $i=2$ }\\
\mathbb{Z}^4 \oplus (\mathbb{Z}/2)^4 & \mbox{ if $i=3$ } \\
\mathbb{Z}^5 \oplus (\mathbb Z/4)^4\oplus (\mathbb{Z}/2)^{14}& 
\mbox{ if $i=4$ }\\
(\mathbb Z/2)^{12}& \mbox{ if $i=5$ }\\
\mathbb{Z} \oplus (\mathbb Z/4)^7\oplus (\mathbb{Z}/2)^{20} & 
\mbox{ if $i=6$ } \\
(\mathbb{Z}/2)^{12} & \mbox{ if $i=2k+1, k>1$ } \\
(\mathbb Z/4)^8\oplus (\mathbb{Z}/2)^{20} & 
\mbox{ if $i=2k, k>3$ }
\end{array} \right.
 \]
We obtain a calculation for the group of gerbes:

$$FGb(T^6/\mathbb{Z}/4)=U(1)^5 \oplus (\mathbb{Z}/2)^4,
\text{  } Gb( T^6/\mathbb{Z}/4)=\mathbb{Z}^4 \oplus
(\mathbb{Z}/2)^4$$

We now offer an interpretation of this calculation. The fact that
$H^1(\mathbb Z/2, M)\cong (\mathbb Z/2)^4$ means that there are
sixteen distinct $M$-equivalence classes of subgroups of order two
in $\Gamma = M\rtimes \mathbb Z/4$. Of these sixteen  classes,
twelve of them have a normalizer of the form $\mathbb Z^2\times\mathbb Z/2$
and four of them have normalizer $\mathbb Z^2\rtimes\mathbb Z/4$, where
$\mathbb Z/4$ acts via multiplication by $-I$ on $\mathbb Z^2$.
In $\Gamma$, the twelve conjugacy classes are identified in pairs,
and so we obtain:

\begin{prop}
The group $\Gamma = M\rtimes\mathbb Z/4$ has ten conjugacy classes of
subgroups of order two. Six of them have normalizers
$\mathbb Z^2\times\mathbb Z/2$ and four of them have
normalizers $\mathbb Z^2\rtimes\mathbb Z/4$.
\end{prop}

Given that $\Gamma$ is a group of finite virtual cohomological dimension
with periodic cohomology,
we can apply a formula due to Brown (see \cite{brown}, page 293)
for the
high--dimensional cohomology of $\Gamma$, expressible as the sum
of the cohomology of the normalizers of elements of order two,
which in this case becomes (for $i>6$):

$$H^i(\Gamma , \mathbb Z)
\cong [H^i(\mathbb Z^2\times\mathbb Z/2, \mathbb Z)]^6
\oplus [H^i(\mathbb Z^2\rtimes\mathbb Z/4, \mathbb Z)]^4
$$
which indeed coincides with our computation.

D. Joyce (\cite{joyce}, pg. 861) has constructed five different desingularizations
of the orbifold $Y_1$. His construction arises from desingularizing the
orbifolds $\mathbb T^2/\{\pm I\}$ associated to the four centralizers
of the form $\mathbb Z^2\rtimes\mathbb Z/4$. Each of them can
be resolved in two possible ways and this leads to the existence
of five different possibilities.
From the point of view of group cohomology,
we see that these centralizers are captured by the high dimensional
cohomology. In fact the Farrell--Tate Cohomology of
the orbifold fundamental group seems to be a natural repository for
the relevant cohomological data. For our example this invariant
is periodic of period two,
and equals
$\widehat H^1 (\Gamma )\cong (\mathbb Z/2)^{12}$,
$\widehat H^2(\Gamma )\cong (\mathbb Z/4 )^8 \oplus (\mathbb Z/2)^{20}$.
It would seem from this calculation that it is naive to expect that
desingularizations can be read off from the cohomology. On the other hand
we can see that the cohomology does contain important information about the
singular set. 

\subsection{The case $Y_2=T^6/(\mathbb{Z}/2)^2$}
Here the action of $(\mathbb{Z}/2)^2$ on $T^6$ is also induced from the
action of $\mathbb{Z}/2^2$  on $\mathbb{C}^3$, given on generators by 
$$\sigma_1(z_1,z_2,z_3)=(-z_1,-z_2,z_3),\;\;\;
\sigma_2(z_1,z_2,z_3)=(-z_1,z_2,-z_3).$$
The induced action of $\mathbb{Z}/2^2$  on $\mathbb{Z}^6$ is given
by the following decomposition of $(\mathbb{Z}/2)^2$ modules:

$$M=(L_1)^2 \oplus (L_2)^2 \oplus (L_3)^2$$
where $L_1$ has rank one and $\sigma_i , i=1,2$, act by multiplication
by $-1$; $L_2$ has rank one, $\sigma_1$ acts by multiplication by
$-1$, $\sigma_2$ acts by multiplication by $1$; and  $L_3$ also has
rank one
and $\sigma_1$ acts by multiplication by $1$ and $\sigma_2$ acts
by multiplication by $-1$.

As before, we need a decomposition of the exterior powers of $M$;
however as it breaks up into one dimensional summands this is
easy to obtain. Note that $L_i\otimes L_i \cong \mathbb Z$,
whereas $L_i\otimes L_j\cong L_k$, where $\{i,j,k\} = \{1,2,3\}$.

\[H^i(\mathbb{Z}^6,\mathbb{Z})= \left\{ \begin{array}{ll}
\mathbb{Z} & \mbox{ if $i=0,6$ }\\ L_1^2 \oplus L_2^2 \oplus
L_3^2 &
\mbox{ if $i=1,5$ } \\
\mathbb{Z}^3 \oplus L_1^4 \oplus L_2^4
\oplus L_3^4 & \mbox{ if $i=2,4$ }\\
\mathbb{Z}^8 \oplus L_1^4 \oplus L_2^4 \oplus L_3^4
 & \mbox{ if $i=3$}
\\ 0 & \mbox{ if $i>6$}
\end{array} \right.
 \]
We now compute the first few cohomology groups for $\mathbb Z^6\rtimes
(\mathbb Z/2)^2$. First we observe that if $G=(\mathbb Z/2)^2$, then
for each $i=1,2,3$ there exists an index two subgroup $H_i\subset G$
such that $L_i$ fits into a short exact sequence
$0\to L_i\to \mathbb Z [G/H_i]\to\mathbb Z\to 0$. From this we conclude
that $H^1(G, L_i)\cong\mathbb Z/2$ and  $H^2(G, L_i)\cong \mathbb Z/2$.
Hence we have that

\[ H^i(\mathbb{Z}^6 \rtimes \mathbb{Z}/2^2,\mathbb{Z})= \left\{ \begin{array}{ll} \mathbb{Z} &
\mbox{ if $i=0$}\\ 0 & \mbox{ if $i=1$} \\ \mathbb{Z}^3 \oplus
(\mathbb{Z}/2)^{6} & \mbox{ if $i=2$} \\
\mathbb{Z}^8 \oplus (\mathbb{Z}/2)^{18} & \mbox{ if $i=3$} 
\end{array} \right.
 \]
From this we derive the desired information about gerbes:
$$FGb(T^6/(\mathbb{Z}/2)^2)=U(1)^3 \oplus (\mathbb{Z}/2)^{18},
\text{  } Gb(( T^6/(\mathbb{Z}/2)^2)=\mathbb{Z}^8 \oplus
(\mathbb{Z}/2)^{18}$$

Of course we could write down the entire cohomology with integral
coefficients, but that would not be particularly illuminating.
However we can use $\mathbb F_2$ coefficients instead; the spectral
sequence still collapses and the complete ring structure 
is given up to filtration by

$$H^*(\mathbb Z^6\rtimes (\mathbb Z/2)^2, \mathbb F_2)
\cong \mathbb F_2 [u, v]\otimes \wedge (x_1, x_2, x_3, x_4, x_5, x_6)$$
where the generators are all in degree one. Notice that the fixed point
set of this action consists of $2^6 = 64$ isolated points, and that
the total dimension of the cohomology of the torus is also $2^6$. This
is of course consistent with the collapse of the spectral sequence and
the localization theorem.

This can be stated as a special case of a general theorem:

\begin{thm}
Suppose that $G$ is a finite group and we are
given a homomorphism 
$$G\to GL_1(\mathbb Z)^n\subset GL_n(\mathbb Z).$$
If $\Gamma = M\rtimes G$ where $M$ is the resulting $\mathbb Z [G]$
lattice, then
$$H^*(\Gamma, \mathbb F_2)\cong H^*(G,\mathbb F_2)\otimes\wedge^*(M).$$
\end{thm}

From this example we see that the cohomology will grow very rapidly, as
to be expected in a non--periodic situation. This is reflected in the
large number of desingularizations which have been described in
\cite{joyce}. Again we see that the cohomological computation provides
insight as to the complicated possibilities for desingularization, but
their precise determination is quite complex.

\end{document}